\documentclass[12pt,twoside,a4paper]{amsart}
\usepackage{amssymb}
\date{\today}

\def\deg{\text{deg}\,}

\def\w{\wedge}

\def\dbar{\bar\partial}

\def\R{{\mathbb R}}
\def\C{{\mathbb C}}
\def\P{{\mathbb P}}
\def\T{{\mathbb T}}

\def\G{{\mathcal G}}

\def\ot{\leftarrow}

\def\D{{\mathcal D}}
\def\S{{\mathcal S}}

\def\supp{{\rm supp\, }}

\def\Im{{\rm Im\, }}

\def\Ker{{\rm Ker\,  }}
\def\Dom{{\rm Dom\,  }}

\def\E{{\mathcal E}}

\def\O{{\mathcal O}}
\def\L{{\mathcal L}}

\def\Re{{\rm Re\,  }}

\def\L{{\mathcal L}}

\def\be{\begin{equation}}
\def\ee{\end{equation}}

\def\A{{\mathcal A}}
\def\sr{\stackrel}

\newtheorem{thm}{Theorem}[section]
\newtheorem{lma}[thm]{Lemma}
\newtheorem{cor}[thm]{Corollary}
\newtheorem{prop}[thm]{Proposition}

\theoremstyle{definition}

\newtheorem{df}{Definition}

\theoremstyle{remark}

\newtheorem{preremark}{Remark}
\newtheorem{preex}{Example}

\newenvironment{remark}{\begin{preremark}}{\qed\end{preremark}}
\newenvironment{ex}{\begin{preex}}{\qed\end{preex}}

\numberwithin{equation}{section}

\begin{document}

\title[Operators with smooth functional \ldots]
{Operators with  smooth
functional calculi}

\date{\today}

\author{Mats Andersson \& H\aa kan Samuelsson \& Sebastian Sandberg}

\address{Department of Mathematics\\Chalmers University of Technology 
and the University of G\"oteborg\\S-412 96 G\"OTEBORG\\SWEDEN}

\email{matsa@math.chalmers.se}


\begin{abstract}
We introduce  a class of (tuples of commuting) unbounded operators on a Banach space,
admitting  smooth functional calculi, that  contains all  operators
of Helffer-Sj\"ostrand type  and is closed under the action of  smooth proper
mappings. Moreover, the class  is closed under tensor product of commuting operators.
In general  an  operator in this class  has no resolvent in the
usual sense so the spectrum must be defined in terms of the functional calculus.
We also consider invariant subspaces and spectral decompositions.
\end{abstract}

\maketitle

\section{Introduction}


In this paper we study unbounded operators on a Banach space $X$ that
admit smooth functional calculi, although they do  not necessarily
have resolvents. 
Throughout this paper $X$ is a complex Banach space, 
$\L(X)$ is  the space of bounded linear operators
on $X$, and $e_X$ denotes the identity operator.

Let  $a$ be  a closed (densely defined) operator 
with real spectrum and with the property  that for each compact set $K\subset\subset\C$ 
there are  $N_K$ and $C_K$ such that 
\begin{equation}\label{tempupp}
\|\omega_{z-a}\|\le C_K|\Im z|^{-N_K},\quad z\in K\setminus\R,
\end{equation}
where $\omega_{z-a}$ is the resolvent form
$
\omega_{z-a}=(z-a)^{-1}dz/2\pi i.
$
Then there is  
a continuous multiplicative mapping
$
[a]\colon \D(\R)\to\L(X), 
$
defined by 
\begin{equation}\label{rep}
[a](\phi)=\int \omega_{z-a}\w\dbar\tilde\phi,
\end{equation}
where  $\tilde\phi$ is an almost holomorphic extension to
$\C$ of $\phi$ with compact support.
This was  done by Dynkin, \cite{D},  for bounded operators  $a$ and
for unbounded operators by Helffer and Sj\"ostrand, \cite{HeSj}.
If $a$ is bounded,  $[a]$ acts on all smooth functions $\phi$ on $\R$ and it
coincides with the  holomorphic functional calculus if $\phi$ is holomorphic
in a neighborhood of the spectrum. 
In general,  $[a]$ has a continuous extension to  the algebra $\G$ of all
smooth functions on $\R$ that are holomorphic at infinity,
in particular, to each $r_z(\xi)=1/(z-\xi)$ for $z\in\C\setminus\R$,
and $[a](r_z)=(z-a)^{-1}$.
Conversely, it was proved in \cite{AS} that if there exists  such a multiplicative mapping  $[a]$
such that, in addition,  
\begin{equation}\label{buss}
\cup_{\phi\in\D(\R)} \Im [a](\phi)\ {\rm is\  dense},\quad  \cap _{\phi\in\D(\R)}\Ker [a](\phi)=\{0\},
\end{equation}
and $[a]$ extends continuously to $\G$, 
then there is a closed operator $a$ satisfying \eqref{tempupp} and such that
\eqref{rep} holds, see Theorem~\ref{gsats}  for the precise statement.
However, in many cases there exists such a smooth  functional
calculus although the resolvent does not exist at all.
For example, let   $a$ be   multiplication with
$\xi\mapsto \xi(2+\sin\xi^3)$  on $X=H^1(\R)$.
Then the resolvent set  is empty, but nevertheless $a$ admits a smooth 
functional calculus $\D(\R)\to\L(X)$, and \eqref{buss} holds.

We take the existence of a smooth functional calculus as our starting point,
and introduce the notion of a {\it hyperoperator}, (with respect to smooth functions).
It is  a multiplicative  $\L(X)$-valued distribution  $A$ on $\R$
such that \eqref{buss} holds. 
This additional requirement means that  $A(1)=e_X$ in a  weak sense.
The spectrum of $A$ is defined as the support of
the distribution.
A closable operator (tuple of commuting closable operators) 
defined on a dense subspace $D$ is a {\it weak hyperoperator}, who, if
 $a$ admits an $\E$ functional calculus with respect to $D$,
i.e., a multiplicative continuous mapping $\E(\R^n)\to \L(D)$,
where $\L(D)$ is the set of closable operators mappings $D\to D$.
Roughly speaking this means that
each $x\in D$ has real and compact local spectrum with respect to $D$.
If $a$ is a who and $f$ is any smooth mapping then $f(a)$ is again  a who. 
It turns out that for any hyperoperator $A$  there is an associated  who $a$.
If $f$ is proper, then  the push-forward $B=f_*A$ of $A$ is a hyperoperator
and $b=f(a)$ is the  who associated to $B$. Conversely, a who
 $a$ is (or corresponds to) a hyperoperator if and only if for
each $\phi\in\D(\R^n)$, $\phi(a)$ extends to a bounded operator on $X$.
Moreover, $a$ is bounded (extends  to a bounded operator) if and only if for each
 $f\in\E(\R^n)$, $f(a)$ extends to a bounded operator on $X$. 

It is a well-known problem to find a suitable  definition of commutativity
for unbounded operators to get a reasonable theory.
We will  consider hyperoperators on $\R^n$ as well, 
with a completely analogous definition. 
For instance, if  $A_1$ and $A_2$ are hyperoperators in $\R$, with associated  whos
$a_1$ and $a_2$, 
commuting in the functional 
calculus sense, then $A=A_1\otimes A_2$ is a new hyperoperator in $\R^2$,
and $a=(a_1,a_2)$ is the associated  who.    However, it is not true that
each hyperoperator in $\R^2$ appears  in this way.
Similar phenomena  hold  for the unbounded analogs of a 
 commuting tuple of bounded operators that are studied in
e.g.,  \cite{IV}, \cite{Sam}, \cite{Vas}, and \cite{Vas1}. 
This gives support for the idea that  a reasonable notion of
``commuting tuple of unbounded operators'' 
must be considered as an object in  its own.
Weaker forms of commutativity of unbounded operators are studied
in \cite{Sch1}, \cite{Sch2}, \cite{Sch3}, and \cite{SchFri}.

One can think of \eqref{rep} as meaning that 
\begin{equation}\label{ekv}
\dbar\omega_{z-a}=[a],
\end{equation}
where $[a]$ is the operator-valued distribution $\phi\mapsto \phi(a)$.
For a general hyperoperator the resolvent form does not exist, but we present
other solutions to \eqref{ekv}  such that   representations
like \eqref{rep} still holds.

\tableofcontents

\section{Notation and some preliminaries}\label{preliminaries}

Any  closed  (densely defined) operator $a$ on  $X$, 
has a well-defined resolvent set $\rho(a)$ which is an open (possibly empty)
subset of the extended plane $\widehat\C$.
The spectrum
of $a$ is the set $\sigma(a)=\widehat\C\setminus\rho(a)$.
Moreover, the operator $a$ is bounded if and only if its spectrum
is contained in $\C$.
For any automorphism $\phi(\zeta)$ of  $\widehat\C$ such that
$\phi^{-1}(\infty)$ is not in the point spectrum of $a$, $\phi(a)$ is a
well-defined closed operator,  and the spectral mapping property
$
\phi(\sigma(a))=\sigma(\phi(a))
$
holds. The automorphism 
\begin{equation}\label{cayley}
C(\zeta)=\frac{\zeta+i}{\zeta-i}, \quad C^{-1}(\tau)=i\frac{\tau+1}{\tau-1},
\end{equation}
maps $\widehat\R$ bijectively onto to the unit circle $\T$. It induces
the Cayley transform which establishes a one-to-one correspondence between
closed  operators with spectrum contained in $\R$, and
bounded operators $b$ with spectrum contained in $\T$ such that 
$b-e_X$ is injective.

\smallskip
If $a$ is a densely defined operator on $X$, then it is closable if
there is a closed operator $a'$ such that $a\subset a'$, i.e., that
the graph of $a$ is contained in the graph of $a'$. 
In that case the closure of the graph of $a$ is the graph of a
(closed) operator called the closure $\bar a$ of $a$. If
$a$ has a bounded extension, then it is equal to $\bar a$.
\smallskip

We let $H^k(\R^n)$ denote the Sobolev space consisting of
all functions in $L^2(\R^n)$ such that all derivatives
up to order $k$ belongs to $L^2(\R^n)$ as well.

\subsection{The Dynkin-Helffer-Sj\"ostrand functional calculus}
For any  $\phi\in\D(\R)$ one can find an extension $\tilde\phi$ to $\C$ such that
$$
\dbar\tilde\phi(\zeta)= \O(|\Im\zeta|^\infty);
$$
such a $\tilde\phi$ is called an almost holomorphic extension of $\phi$.
Moreover, if $\tilde K$ is a complex neighborhood of $\supp\phi$, one may assume 
that $\tilde\phi$ has support in $\tilde K$.
Now let  $a$ be a closed operator with real spectrum such that
\eqref{tempupp}  holds, such an operator will be referred to as an HS operator.
Then clearly the integral in \eqref{rep} converges,  and it turns out
to be independent of the choice of almost holomorphic extension. 
The multiplicativity follows from an application of the  resolvent identity
$$
\frac{1}{w-a}\frac{1}{z-a}=\frac{1}{z-w}\frac{1}{w-a}+\frac{1}{w-z}\frac{1}{z-a}.
$$
It is easy to see  that $[a]$ is continuous in the sense that
$[a](\phi_j)\to 0$ in operator norm if $\phi_j\to 0$ in $\D(\R)$.
It also follows that the support of $[a]$ coincides with $\sigma(a)\cap\C$.
Moreover, we claim that 
\begin{equation}\label{buster}
[a](\xi\phi)x=[a](\phi)ax,\quad   x\in\Dom(a).
\end{equation}
This is of course well-known, but for further reference we sketch a proof.
From the resolvent identity we have, assuming that $-i$ is outside the 
support of $\tilde\phi$, 
\begin{multline*}
\int\Big(\frac{1}{a+i}-\frac{1}{z+i}\Big)\frac{dz}{z-a}\w\dbar\tilde\phi(z)=\\
\int\frac{dz}{(a+i)(z+i)}\w\dbar\tilde\phi(z)=
-\int\dbar
\Big(\frac{dz}{(a+i)(z+i)}\w\tilde\phi(z)\Big)=0,
\end{multline*}
where the last equality follows from Stokes' theorem. Thus we have
\begin{equation}\label{buster1}
\frac{1}{a+i}[a](\phi)=[a](\phi)\frac{1}{a+i}=[a](\frac{1}{\xi+i}\phi(\xi)).
\end{equation}
Replacing $\phi(\xi)$ by $(\xi+i)\phi(\xi)$ we get
\begin{equation}\label{buster2}
\frac{1}{a+i}[a]((\xi+i)\phi)=  [a]((\xi+i)\phi) \frac{1}{a+i}=   [a](\phi).
\end{equation}
If $x\in\Dom(a)$ we therefore have
$[a]((\xi+i)\phi)x=[a](\phi)(a+i)x$, which implies \eqref{buster}.

\begin{ex}\label{volterra}
Let $a$ be a closed operator with spectrum equal to
$\{\infty\}$. For instance one can take 
 the inverse of the Volterra operator. 
Then clearly \eqref{tempupp} holds,
but the resulting  multiplicative mapping  $[a]$   is identically $0$.
\end{ex}

If $a_1,\ldots,a_n$ is a tuple of HS operators such that their resolvents
(anti-) commute, i.e.,
$
\omega_{\zeta_j-a_j}\wedge\omega_{\zeta_k-a_k}=-\omega_{\zeta_k-a_k}\wedge
\omega_{\zeta_j-a_j},$
for  $\zeta_j,\zeta_k\in\C\setminus\R$, then
$[a]=[a_1]\otimes [a_2]\cdots\otimes [a_n]\in\D'(\R^n,\L(X))$ is multiplicative.
This follows by  simple abstract considerations, but it can also be realized
explicitly as
$$
[a](\phi)=\int
\omega_{\zeta_1-a_1}\wedge\ldots\w\omega_{\zeta_n-a_n}\w\dbar_{\zeta_n}\cdots\dbar_{\zeta_1}\tilde\phi(\zeta),
$$
where $\tilde\phi$ is a special almost holomorphic extension to
$\C^n$ with compact support  as in  \cite{AS}, i.e., such that
\begin{equation}\label{asspecial}
\dbar_{\zeta_1}\cdots\dbar_{\zeta_m}\tilde\phi(\zeta)=\O(|\Im\zeta_1|^\infty\cdots|\Im\zeta_m|^\infty).
\end{equation}

\subsection{Commuting bounded operators}

Let $a=(a_1,\ldots,a_n)$  be a commuting tuple of bounded operators on $X$.
If the Taylor spectrum $\sigma(a)$ is contained in $\R^n$, then it coincides
with the spectrum of $a$ with respect to the commutative Banach algebra
$(a)$ generated by $a$.
If the tuple $a$ has real spectrum, then  we say that $a$ admits a 
smooth functional calculus if the real-analytic functional
calculus  $[a]\colon \O(\R^n)\to\L(X)$
has a continuous extension to a mapping $[a]\colon\E(\R^n)\to\L(X)$.
Since $C^\omega(\R^n)$ is dense in $\E(\R^n)$, the extension
is then unique and multiplicative, and in fact it extends to
$\E(\sigma(a))\to\L(X)$.
 The existence of such an extension
is equivalent to that
$\exp(ia t)$ has polynomial growth in $t\in\R^n$, see, e.g., \cite{A2};  it is 
also equivalent to that the resolvent satisfies
$$
\|\omega_{z-a}\|\le C |\Im z|^{-M},
$$
for some $M>0$.

\smallskip
If $a$  has  non-real (Taylor) spectrum $\sigma(a)$, then there is in general
no unique extension of the holomorphic functional calculus. 
For instance, let $b$ be a nilpotent operator
and let $A(\phi)= \tilde\phi(b,0)$  and $B(\phi)=\tilde\phi(b,b)$ respectively,
where $\tilde\phi(z,\bar z)= \phi(z)$ for real-analytic $\phi$
(only a finite Taylor expansion is needed). Then
$A$ and $B$ extend to two different multiplicative mappings 
$\E(\C)\to\L(X)$ which both extend the holomorphic functional
calculus.
In general, a  possible smooth functional calculus is uniquely determined by
the image of $\bar{z}$ (or $\bar{z}_j$ if we have an $n$-tuple
of commuting operators). In our situation the bounded (tuples of) operators that appear
are like  $b=A(\phi)$ for a possibly complex-valued $\phi$, and then
we have a natural conjugated operator, namely $b^*=A(\bar\phi)$.
A smooth functional calculus for such an operator $b$ is then understood
to map $\bar{z}$ to $b^*$.
If $f(z)=\tilde f(z,\bar z)=\hat f(\Re z,\Im z)$, then
$f(b)=\tilde f(b,b^*)=\hat f((b+b^*)/2,(b-b^*)/2i)$,
and therefore we can reduce to the case of real-valued functions $\phi$.
\smallskip

We conclude this section with the following useful observation.
\begin{lma}\label{hololemma}
If $A$ is a linear and multiplicative mapping $\mathcal{D}(\R)\rightarrow \mathcal{L}(X)$
then, for any $\chi\in \D(\R)$, $z \mapsto A(\chi(\xi)/(z-\xi ))$ is strongly holomorphic in 
$\mathbb{C}\setminus \mathbb{R}$.
\end{lma}
\begin{proof}
Let $\tilde{\chi}\in \D(\R)$ be identically $1$ on $\mbox{supp}\, \chi$.
From linearity and multiplicativity we get
\begin{equation}\label{eq2}
A(\frac{\chi(\xi)}{z-\xi })=A(\frac{\chi(\xi )}{z_0-\xi })-(z-z_0)A(\frac{\chi(\xi )}{z-\xi })
A(\frac{\tilde\chi(\xi)}{z_0-\xi }).
\end{equation}
Letting $\|A(\chi(\xi )/(z_0-\xi ))\|=C$ and  $\|A(\tilde\chi(\xi )/(z_0-\xi ))\|=\tilde C$
we see that 
\[
\|A(\chi(\xi )/(z-\xi ))\|\leq C+|z-z_0|\tilde C\|A(\chi(\xi )/(z-\xi ))\|
\]
and so 
$\|A(\chi(\xi )/(z-\xi ))\|\leq C/(1-|z-z_0|\tilde C).$
Thus $\|A(\chi(\xi )/(z-\xi ))\|$ is locally uniformly bounded in $z$. From \eqref{eq2}
it now follows that $A(\chi(\xi)/(z-\xi))$ is strongly 
continuous at $z_0$. With this fact in
mind it follows immediately from \eqref{eq2} that 
\[
\frac{1}{z-z_0}(A(\frac{\chi(\xi )}{z-\xi })-A(\frac{\chi(\xi )}{z_0-\xi }))\rightarrow
- A(\frac{\chi(\xi )}{(z_0-\xi )^2}), \,\,\, z\rightarrow z_0,
\]
in operator norm.
\end{proof}

\section{Definition and basic properties}

We say that a linear mapping 
$A\colon \D(\R^n)\rightarrow \L(X)$ is continuous,
$A\in \mathcal{D}'(\R^n,\L(X))$, if
$A(\phi_j)\to 0$ in operator norm when  $\phi_j\to 0$ in $\D(\R^n)$.
As for ordinary distributions
it follows immediately that $A$ has finite order on compact subsets, i.e.,
 for any compact $K \subset \R^n$ there is a constant  $C_{K}$ and 
a non-negative integer $M_K$ such that 
$$
\| A(\phi) \| \leq C_K\sum_{|\alpha|\leq M_K}\sup_K|\partial^{\alpha}\phi|
$$
for all $\phi \in \D(\R^n)$ with support in $K$.

\begin{df}\label{hodef}
A continuous multiplicative mapping  $A\colon \D(\R^n)\to \L(X)$ is a hyperoperator on $\R^n$,
$A\in H_{\D(\R^n)}(X)$, if
$$
(i)\quad D_A=\cup \Im A(\phi)\ {\rm is\ dense\ in}\  X, \quad\quad {\rm and}
$$
$$
(ii)\quad N=\cap\Ker A(\phi)=\{0\}.
$$
\end{df}

If $a$ is an  HS operator such that  $[a]$ satisfies $(i)$ and $(ii)$,  then 
$[a]$ is a hyperoperator.    
If $a$ is bounded  (or a commuting tuple of bounded operators), then
$[a](\phi)=\phi(a)$.
It is readily checked that  the operator (tuple of operators) $0_X$ gives rise to the
hyperoperator   $[0_X]$, defined  by  $[0_X](\phi)=\phi(0)e_X$.
In the same way, $[e_X](\phi)=\phi(1)e_X$.

\begin{remark} 
Let  $A\colon \D(\R^n)\to \L(X)$ be  a continuous multiplicative mapping.
If $A$ has  compact support, i.e.,   $A$ has a  continuous
extension to $\E(\R^n)$, then   $(i)$ and $(ii)$ hold if and only if $A(1)=e_X$. 
In fact,
let $\chi_N$ be a sequence in $\D(\R^n)$ that tends to $1$ in
$\E(\R^n)$. If now $A(1)=e_X$, then for any $x\in X$ we have
that $x=A(1)x=\lim A(\chi_N)x$,  and hence $(i)$ holds.
In the same way, if $A(\chi_N)x=0$ for all $N$, then
$x=A(1)x=0$ so that $(ii)$ holds as well. Conversely,
if $x\in D_A$, then $x=A(\phi)z$ and therefore
$A(1)x=A(1)A(\phi)z=A(1\cdot\phi)z=A(\phi)z=x$.
If $D_A$ is dense it follows that $A(1)=e_X$.
Therefore it is natural to think of $(i)$ and $(ii)$ as
a weak form of saying that $A(1)=e_X$.
\end{remark}

We say that $\chi_N\in\D(\R^n)$ is an exhausting sequence if 
$0\le\chi_N\le 1$, $\chi_N\nearrow 1$,  and the compact sets
$K_N=\{\chi_N=1\}$ form an exhausting sequence of compact sets;
i.e., $K_N\subset int(K_{N+1})$ and $\cup K_N=\R^n$.

\begin{lma}\label{lma1}
Suppose that $\chi_N$ is an exhausting sequence in $\R^n$ and
$A\in H_{\D(\R^n)}(X)$. Then
$
\cup\Im A(\chi_N)=D_A.
$
\end{lma}

\begin{proof}
If $\phi\in\D(\R^n)$, then $\chi_N\phi=\phi$ if $N$ is large enough, and therefore
$$
A(\chi_N)A(\phi)=A(\chi_N\phi)=A(\phi),
$$
which shows that $A(\chi_N)$ is the identity on $\Im A(\phi)$.
Thus $\Im A(\chi_N)\supset \Im A(\phi)$.
\end{proof}

\begin{prop}\label{tensorprop}
Assume that $A_1$ and $A_2$ are hyperoperators in
$\R^n$ and $\R^m$,  respectively, and that they are commuting,
i.e.,
$$
A_1(\phi)A_2(\psi)=A_2(\psi)A_1(\phi), \quad \phi\in\D(\R^n),\ \psi \in\D(\R^m).
$$
Then 
$A=A_1\otimes A_2$ 
is a hyperoperator in $\R^{n+m}$ and
$
D_A=D_{A_1}\cap D_{A_2}.
$
\end{prop}

In particular it follows that $D_{A_1}\cap D_{A_2}$
is dense as soon as $A_1$ and $A_2$ are commuting.

\begin{proof}
The tensor product $A$ is defined as usual for distributions; thus 
$A(\phi\otimes\psi)=A_1(\phi) A_2(\psi)$,  and it is  extended
to $\D(\R^{n+m})$ by linearity and  continuity.
The assumption on commutativity implies that $A$ is multiplicative. 
If $0=A(\phi\otimes\psi)x=A_1(\phi) A_2(\psi)x$ for all $\phi$ and $\psi$
it follows from condition $(ii)$ for $A_1$ and $A_2$ that $x=0$.
Thus $(ii)$ holds for $A$.
Given $x\in X$ we can find $y$ and $\phi$ such that $\|x-A_1(\phi)y\|<\epsilon/2$.
In the same way we can find $z$ and $\psi$ such that
$\|y-A_2(\psi)z\|<\epsilon/(2\|A_1(\phi)\|)$. It follows that
$\|x-A(\phi\otimes\psi)z\|<\epsilon$. Thus $D_A$ is dense in $X$.
On the other hand, since $\chi_N\otimes\chi'_M$ is an exhausting
sequence in $\R^{n+m}$  if $\chi_N$ and $\chi'_M$ are exhausting
sequences in $\R^n$ and $\R^m$, respectively, it follows that 
$x\in D_A$ if and only if  $A(\chi_N\otimes\chi_M)x=x$ for sufficiently  large $N$ and $M$,
and this in turn holds if and only if $x\in D_{A_1}\cap D_{A_2}$.
\end{proof}

If a hyperoperator $A$ in $\R^{n+m}$ is the tensor product $A_1\otimes A_2$ 
of two commuting, multiplicative $\mathcal{L}(X)$-valued distributions in 
$\R^n$ and $\R^m$, then each $A_j$ is indeed a hyperoperator. 
In fact,  since $\chi_N\otimes\chi'_M$ is exhausting in $\R^{n+m}$,
$\cup \mbox{Im}\,A(\chi_N\otimes\chi'_M)=\cup \mbox{Im}\,A_1(\chi_N)A_2(\chi'_M)=
\cup \mbox{Im}\,A_2(\chi'_M)A_1(\chi_N)$ is dense, so   $A_j$ satisfy condition $(i)$.
If  $A_j(\phi)x=0$ for all $\phi \in \mathcal{D}(\R^n)$,  then 
$A(\phi\otimes\phi')x=0$ for all $\phi, \phi'$. Therefore
$A(\psi)x=0$ for all $\psi\in\D(\R^{n+m})$, so   $x=0$. Hence $A_j$
satisfies $(ii)$.

\begin{prop}
If $A\in H_{\D(\R^n)}(X)$ and $f\in\E(\R^n,\R^m)$ is a proper mapping, then
the push-forward $B=f_*A\in\D'(\R^m,\L(X))$ is a hyperoperator,
and $D_B=D_A$.
\end{prop}

\begin{proof}
Since $f$ is proper, $f^*\colon \D(\R^m)\to\D(\R^n)$ and hence
$f_*A$, defined by $f_*A(\phi)=A(f^*\phi)=A(\phi\circ f)$, is
a multiplicative distribution. 
If $\chi_N$ is an exhausting sequence in $\R^m$, since $f$ is proper,  then 
$\chi_N\circ f$ is an exhausting sequence in $\R^n$. Therefore, 
$$
D_B=\cup_N\Im f_*A(\chi_N)=\cup_N \Im A(\chi_N\circ f)=D_A
$$
according to Lemma~\ref{lma1}. Thus $f_*A$ satisfies $(i)$. 
Finally, suppose that $f_*A(\psi)y=0$ for all $\psi\in\D(\R^m)$.
For fixed $\phi\in\D(\R^n)$ and large $N$, then
$$
A(\phi)y=A\big(\phi(\chi_N\circ f)\big)y=A(\phi)A(\chi_N\circ f)y=0,
$$
and since $\phi$ is arbitrary, we conclude that $y=0$. Thus
$f_*A$ is a hyperoperator.
\end{proof}

It is easy to check  that any hyperoperator $A$ extends to a multiplicative
mapping on the algebra $\D(\R^n)\oplus\C$ of smooth
functions that are constant outside some compact set, just by letting
$A(h)=h(\infty)e_X + A(h-h(\infty))$.
If $\phi$ has compact support, then $h=f\circ\phi$ is in this algebra, and therefore
we have

\begin{prop}\label{sallad}
Assume that $A\in H_{\mathcal{D}(\mathbb{R}^n)}(X)$ and $\phi\in\mathcal{D}
(\mathbb{R}^n,\mathbb{R}^m)$. Then
the bounded operator $\phi(a)=A(\phi)$ admits a $\mathcal{E}$-functional calculus
that extends the holomorphic (real-analytic)  functional calculus, defined by
$f\mapsto f(0)e_X+ A(f\circ\phi-f(0))$.
\end{prop}

\section{Weak hyperoperators}

We shall now see that for each hyperoperator $A$ there is an associated
closable operator $a$ on $D_A$. We will use the operator $a$ to 
model the definition of a {\em weak hyperoperator}, see Definition \ref{who} below.



\smallskip

Let $A$ be a hyperoperator in $\R^n$ and let $f\colon\R^n\to\R^m$ be any smooth mapping.
If $x\in D_A$ and $x=A(\phi)y$ we define
$f(a)x=A(f\phi)y$. If $\chi=1$ in a neighborhood  of $\supp \phi$, then
$f(a)x=A(f\chi\phi)y=A(f\chi)A(\phi)y=A(f\chi)x$; thus
$f(a)x=A(f\chi)x$ and  in particular $f(a)$ is a well-defined densely
defined operator.
Also observe that if  $\phi\in\D(\R^n)$, then
$\phi(a)x=A(\phi)x$ for all $x\in D_A$.

\smallskip

For any $x\in X$ we let $\sigma_x(A)$ be the support of the $X$-valued
distribution $\phi\mapsto A(\phi)x$; this is the local spectrum at $x$. 
If $K\subset\R^n$ is compact, we let 
$$
D_{A,K}=\{x\in X;\  \sigma_x(A)\subset K\}.
$$
It is readily checked that
$D_A=\cup_K D_{A,K}$.

\begin{prop}\label{tatt}
Assume that $A$ is  a hyperoperator in $\R^n$.

\smallskip
\noindent (a)\quad  If $f\in\E(\R^n,\R)$, then  $f(a)$ maps $D_A\to D_A$, and 
if $g\in\E(\R^n,\R)$, then $g(a)f(a)=(fg)(a)$ on $D_A$.

\smallskip

\noindent (b)\quad If $f\in\E(\R^n,\R^m)$,  then
$f(a)$ is a closable operator (tuple of operators).

\smallskip
\noindent (c) If $f_k\to f$ in $\E(\R^m,\R^p)$, then
$f_k(a)x\to f(a)x$ for all $x\in D_A$.

\smallskip
\noindent (d) If $x_j\in D_{A,K}$ for some fixed compact set $K$
and $x_j\to x$  in $X$, then $f(a) x_j\to f(a) x$.
\end{prop}


\begin{proof}
If $x\in D_A$ and $y=f(a)x$, then $y=A(\chi f)x$ for an appropriate
$\chi$ and hence by definition $y\in D_A$. Moreover,
\begin{multline*}
g(a)f(a)x=g(a)A(\chi f)x=A(\tilde{\chi} g)A(\chi f)x=A(\tilde{\chi}\chi fg)x=\\
A(\chi fg)x=(fg)(a)x,
\end{multline*}
if $A(\chi)x=x$ and $\tilde{\chi}=1$ on the support of $\chi$. Thus (a) holds.

We can always take the closure of the graph of $f(a)$ in $X^n\times X^m$.
If $x_k\in D_A$,  $x_k\to x$,  and $f(a)x_k\to y$, then
for any $\psi\in\D(\R^n)$,
$A(\psi)f(a)x_k\to  A(\psi)y$; but also
$A(\psi)f(a)x_k=A(\psi f)x_k\to A(\psi f)x$,
so $A(\psi)y=A(\psi f)x$.
Because of condition $(ii)$ we have that
$y$ is then uniquely determined by $x$, and
hence the closure is a graph.
Thus $(b)$ is proved.

Given $x\in D_A$, take $\chi$ such that 
$A(\chi)x=x$. Since  $f_k\chi\to f\chi$ in $\D(\R^n)$, we have that
$f_k(a)x=A(f_k\chi)x\to A(\chi f)x=f(a)x$. Thus  $(c)$ holds.
For the last statement, observe that 
$A(\chi)x_k=x_k$ if $\chi=1$ in a neighborhood of $K$. Hence 
$f(a)x_k=A(f\chi)x_k\to A(f\chi)x= f(a)x$.
\end{proof}

Notice that if  $\phi\in\D(\R^n)$, then the closure of $\phi(a)$ is equal to the
bounded operator $A(\phi)$. 
Moreover,  the closure of $1(a)$ is equal to $e_X$ and the closure of
$0(a)$ is equal to $0_X$.
Applying Proposition~\ref{tatt}  to the  mapping $f(\xi)=\xi$, we find that
$a$ has a meaning as a densely defined closable operator (tuple of
closable operators $a_j=f_j(a)$, where $f_j(\xi)=\xi_j$,  that commute
on $D_A$).

In view of this proposition it is natural to introduce a special 
class of densely defined linear operators. If $D$ is
a dense subspace, let $\L(D)$ be the set of closable
linear operators $D\to D$.

\begin{df}\label{who}
Let $c=(c,D)$ be a linear operator mapping 
the dense subspace $D$ of $X$ into itself.
Moreover, assume that $c$ is closable, and that
there is a linear and  multiplicative mapping 
$\E(\R^n)\to\L(D)$, 
that extends the trivial one on  polynomials,  and such that
$h_k(c)x\to h(c)x$, for $x\in D$ if $h_k\to h$ in $\E(\R^n)$. 
Then 
we say that $c$, or rather $(c,D)$,  is a weak hyperoperator, a who.
\end{df}

The sum of two closable operators is not  necessarily closable
(so $\L(D)$ is not a space), 
so part of the requirement is that each polynomial 
$p(c)$ in $c$ is closable.
Moreover, since the  polynomials are dense in $\mathcal{E}$ the
extension to $\E(\R^n)\to  \mathcal{L}(D)$ is unique
if it exists.

Assume that $A$ is a hyperoperator and $f$ is any
smooth mapping.
If $h\in\E(\R^m,\R^p)$ then we can define 
$h(f(a))=(h\circ f)(a)$ on $D_A$.
Therefore $(a,D_A)$ as well as
$(f(a), D_A)$ are whos.
Also notice that if $f$ is proper,  $B=f_* A$, and 
$a$ and $b$ are the associated  whos, then $b=f(a)$.

\smallskip

We say that a who $c=(c,D)$ is (extendable  to) a hyperoperator $A$ if
$D\subset D_A$ and $\phi(c)x=A(\phi)x$ for $x\in D$. If such an $A$ exists
it is unique in view of Proposition~\ref{rattmuff} below. 
In the sequel we will therefore often talk about the hyperoperator $a$, 
meaning that $a$ is the who associated  to some  hyperoperator $A$.

\begin{prop}\label{rattmuff}
Suppose that $A$ and $A'$ are in $H_{\D(\R^n)}(X)$ and  that $D_A\cap D_{A'}$ is dense.
Moreover, assume that there is a dense subspace $D$ of $D_A\cap D_{A'}$ such that
$a_j=a'_j$ on  $D$ and map $D\to D$. Then  $A=A'$.
\end{prop}


\begin{proof} 
If $x\in D$ and $\chi$ is identically $1$ on a large enough set, then
$$
A(\xi_j \chi)x=a_jx=a'_jx=A'(\xi_j\chi)x.
$$
Moreover, if $\tilde\chi$ is $1$  in a neighborhood of $\supp\phi$,
 (recall that $x\in D$ implies that $a_jx\in D$)
\begin{multline*}
A(\xi_k\xi_j\chi)x= A(\xi_k\xi_j\chi\tilde\chi)x=A(\xi_k\chi)A(\xi_j\tilde\chi)x=\\
A(\xi_k\chi) a_jx=a_ka_jx=a'_ka'_jx=\ldots=A'(\xi_k\xi_j\chi)x,
\end{multline*}
and so on, so we get
$A(p\chi)x=A'(p\chi)x$
for all polynomials $p$. If $\phi$ is a test function it follows by the Weierstrass approximation theorem 
that $A(\phi)x=A'(\phi)x$, and hence  $A(\phi)=A'(\phi)$ since $D$ is dense.
\end{proof}

\begin{cor}
If $A$ is a hyperoperator and $f$ is proper, then
$f_*A=[0_X]$ if and only if $f(a)=0$.
In particular, $A=[0_X]$  if and only if $a=0$.
\end{cor}

In fact, if $B=f_*A$, then $b=f(a)$, so $bx=0=0_X x$ for all $x\in D_B=D_A$.
Hence, by the previous proposition,  $B=[0_X]$.

\begin{cor}
If $A,A'$ are commuting  hyperoperators and $a=a'$ on $D_A\cap D_{A'}$, then
$A=A'$.
\end{cor} 

This is just because 
 $D_{A}\cap D_{A'}$ is dense if $A$ and $A'$  commute, cf., Proposition~\ref{tensorprop}.

\smallskip

Assume that $A$ is a hyperoperator and let $h$ be  smooth and constant outside a compact set.
It is easily checked that 
the bounded operator $A(h)=h(\infty)e_X+A(h-h(\infty)$ is the closure of
the densely defined operator $h(a)$. Therefore, cf., Proposition~\ref{sallad},
$$
f(A(\phi))x=(f\circ\phi)(a)x, \quad x\in D_A,
$$
for any smooth $f$ if $\phi$ has compact support.

\begin{prop}\label{sallad3}
Let  $a_0$ be  an HS operator such that 
$[a_0]$  satisfies $(i)$ and $(ii)$ so that $A=[a_0]$ is a hyperoperator.
If $a$ is the associated who, then 
$\bar a =a_0$.
\end{prop}

\begin{proof} Since by assumption $a_0+i$ has a bounded inverse, we have that
$\Dom(a_0)=\Dom(a_0+i)=\Im(a_0+i)^{-1}$. If $x\in D_A$, then 
$x=A(\chi)x$ so by \eqref{buster2}
$$
x=A(\chi)x=\frac{1}{a_0+i}A\big((\xi+i)\chi\big)x,
$$
and hence $x\in\Dom(a_0)$ and by \eqref{buster}, $ax=a_0x$. 
Thus $a\subset a_0$.

Now, if  $x\in\Dom(a_0)$ 
there is some $y$ such that $x=(a_0+i)^{-1}y$. Take $y_k\in D_A\subset\Dom(a_0)$ such that
$y_k\to y$. Then 
$x_k=(a_0+i)^{-1}y_k=(a+i)^{-1}y_k\in D_A$ according to \eqref{buster1}. Thus
$$
ax_k=a_0 x_k=\frac{a_0}{a_0+i}y_k\to \frac{a_0}{a_0+i}y=a_0 x,
$$
since $a_0/(a_0+i)$ is bounded. Therefore, 
$(x,a_0 x)$ belongs to the closure of (the graph of) $a$. 
\end{proof}

It is now easy to  see that there exist  non-trivial hyperoperators.

\begin{ex}\label{ex100}
Let $a$ be the unbounded operator defined as multiplication with $\xi$ on 
$X=H^1(\R_{\xi})$.
It defines a hyperoperator $A=[a]$  and the associated who is $(a,D_A)$ where 
$D_A=\{x\in X;\,\mbox{supp}\,x \subset\subset \mathbb{R}\}$. 
The mapping $f(\xi)=\xi(2+\sin \xi^3)\colon\mathbb{R}\rightarrow \mathbb{R}$
is proper and so $B:=f_*A$ is a hyperoperator with $D_B=D_A$. By definition $B(\phi)$
is multiplication with $\phi \circ f$. The who  $b$ associated to 
$B$ is just multiplication with $f$ because if $x\in D_B$ and $\chi$ is chosen so
that $\mbox{supp}\,x \subset \{\chi \circ f=1\}$ then $bx=B(\xi\chi(\xi))x=
A(f\chi \circ f)x=f(a)x=fx$. 
We claim that $B$ is not $[b_0]$ for any HS operator $b_0$.
If there were such a $b_0$, then  by Proposition~\ref{sallad3},
$\bar b=b_0$ and therefore   there would be a bounded operator $c$ such that
$
c(\bar{b}+i)x=(\bar{b}+i)cx=x
$
 for all $x\in D_B=\mbox{Dom}(b)$. However, then $c$ would have 
to be multiplication with $(f(\xi)+i)^{-1}$ on the image of $D_B$ under $\bar{b}+i$
which again is $D_B$, but this is impossible since multiplication with $(f(\xi)+i)^{-1}$ 
has no bounded extension to all of $X$. 
\end{ex}


\begin{ex}
If $B$ is a hyperoperator in $\R^2$ then the associated  who 
$b$ is equal to $(b_1,b_2)$ where $b_j=\pi_j b$ are whos  as well. 
However it may happen that none of the $b_j$ are
hyperoperators. 
Let  $f_1(\xi)$ be  equal to $\xi$ for $\xi>1$ and $\sin \xi^2$ for $\xi<-1$, and
let $f_2(\xi)=-f_1(-\xi)$. Then   $F=(f_1,f_2)\colon\mathbb{R}\to\mathbb{R}^2$
is proper and therefore  $B:=F_*A$ is a hyperoperator, if 
$A\in H_{\mathcal{D}(\mathbb{R})}(H^1(\mathbb{R}))$ is the hyperoperator that 
sends $\phi$ to multiplication with $\phi$. In this case $b_1=f_1(a)$ and
$b_2=f_2(a)$.
Now,  $\phi(b_j)$ is multiplication with $\phi \circ
f_j$ and this operator has in general no bounded extension to 
$H^1(\mathbb{R})$, so $b_j$ is not a hyperoperator.
 Take for instance 
$\phi \in \mathcal{D}(\mathbb{R})$ such that $\phi'(\xi)=1$ for $-1\leq \xi\leq 1$; then 
$ (\phi \circ f_j)'(\xi)$ is unbounded. 
\end{ex}

%

%

\begin{ex}
Let $(M,\mu)$ be a finite measure space and let $h$ be a real or complex
valued measurable function (tuple of functions) defined a.e.\  with respect
to $\mu$. The operator defined as multiplication with $h$ on $L^p(M,\mu)$, $1\le p<\infty$,
is then a hyperoperator and $\sigma(a)$ (see Section \ref{spect})  
is the essential range of $h$. 
Composing with smooth maps and/or taking tensor
products will not take us outside this class of multiplication operators.
By basic spectral theory any normal operator (tuple of normal commuting
operators) can be viewed as such an operator (tuple of operators) on some
$L^2(M,\mu)$. Therefore, our theory does not add anything to the
usual theory of self-adjoint operators.
\end{ex}



We conclude this section with a result which together with Proposition \ref{tatt} 
characterizes those whos that  are hyperoperators.

\begin{prop}\label{whoho}
Let $a=(a,D)$ be a who  such that
the closure of $\phi(a)$ is bounded on $X$ for all $\phi \in \D(\R^n)$. Assume that 
$\cap_{\phi \in \mathcal{D}}\mbox{Ker}\,\overline{\phi(a)}=\{0\}.$ 
Then the mapping $A$ defined by $A(\phi)=\overline{\phi(a)}$ is a 
hyperoperator with $D_A=\cup_{\phi \in \mathcal{D}} \mbox{Im}\,\overline{\phi(a)}\supseteq D$.
 Moreover if $a'$ is the who 
associated to $A$ then $\overline{a'}=\overline{a}$.
\end{prop}

Let $X$ and $a$ be as in Example \ref{ex100}. 
Then $(a,\D(\R))$ is a who satisfying the hypotheses
of the proposition. The induced hyperoperator is $A=[a]$ and $D_A$ is the space of all $f$
in $X$ with compact support. 

\begin{proof}
We first show that $A$ so defined is a continuous mapping
$\D(\R^n)\to\L(X)$. To this end, we take a compact set $K\subset\R^n$, and a cut-off function
$\chi$ that is $1$ in a neighborhood of $K$. For each $x\in X$ we can define a mapping
$A_x\colon\E(\R^n)\to X$ by
$
A_x f=\overline{(\chi f)(a)}x.
$
For $x\in D_1=D\cap\{|x|\le 1\}$ the mapping $A_x$ is continuous, since
$(a,D)$ is a who. By the Banach-Steinhaus theorem it follows that $\{A_x\}_{x\in D_1}$,
$D_1=D\cap \{|x|\leq 1\}$, is
equi-continuous, which means that 
\begin{equation}\label{hast}
|A_x f|\le C\sum_{|\alpha|\le M}\sup_{K'}|\partial^\alpha f|, 
\end{equation}
for some $C,M$ and $K'$ independent of $x\in D_1$. Applying to $\phi$ with support in
$K$, and using that $D$ is dense, we get
$$
\|\overline{\phi(a)}\|\le  C\sum_{|\alpha|\le M}\sup|\partial^\alpha \phi|.
$$
Thus $A$ is continuous.
The multiplicativity $A(\phi\psi)=A(\phi)A(\psi)$
now follows by continuity, since it holds when applied to  $x\in D$. Moreover, for any 
$x\in D$ the map $\E(\R^n)\ni f \mapsto f(a)x \in X$ is continuous and therefore has compact
support, $\sigma_x(a)$. If $\chi=1$ in a neighborhood of $\sigma_x(a)$ it follows that
$\chi(a)x=1(a)x=x$. Hence $A$ is a hyperoperator with 
$D_A=\cup_{\phi\in \D}\mbox{Im}\,\overline{\phi(a)}\supseteq D$.

It remains to see that $\overline{a'}=\overline{a}$. If $\chi_N$ is  an exhausting sequence,
then 
$\psi_N=\xi\chi_N\rightarrow \xi$ in $\mathcal{E}(\R^n)$ 
and so for $x\in D\subseteq D_A$ we have
\[
ax=\lim_{N\rightarrow \infty}\psi_N(a)x=\lim_{N\rightarrow \infty}A(\psi_N)x=a'x.
\]
Hence $a\subseteq a'$ and so $\overline{a}\subseteq \overline{a'}$. To obtain the converse
inclusion it suffices to show that $\overline{\mbox{Graph}(a)} \supseteq
\mbox{Graph}(a')$. Let $(x,a'x)\in \mbox{Graph}(a')$. Since $x\in D_A$,  there is an 
$N_0$ such that $A(\chi_{N_0})x=x$. Take any sequence $y_j$ in $D$ converging to
$x$ and put $x_j=\chi_{N_0}(a)y_j$. Then $x_j$ is a sequence in $D$ and it also
converges to $x$ since $\chi_{N_0}(a)$ has a bounded extension. It follows that
\[
ax_j=\lim_{N\rightarrow \infty}\psi_N(a)x_j=
\lim_{N\rightarrow \infty}\psi_N(a)\chi_{N_0}(a)y_j=\psi_{N_0}(a)y_j \rightarrow
\overline{\psi_{N_0}(a)}x,
\]
as $j\rightarrow \infty$.
However,  $\overline{\psi_{N_0}(a)}x=a'x$ and hence $(x_j,ax_j)\rightarrow (x,a'x)$, 
that is, $(x,a'x)\in \overline{\mbox{Graph}(a)}$.
\end{proof}


\begin{remark}
Let $(a,D)$ be a who. For each $x\in D$ the mapping $\phi \mapsto \phi(a)x$ is 
a continuous mapping $A_x\colon \E(\R^n)\to X$, and hence it has compact support.
As for a hyperoperator, we can define the local spectrum $\sigma_x(a)$
as this support. If  $D_K=\{x\in D;\ \sigma_x(a)\subset K\}$, then 
clearly $D=\cup D_K$. For each $x\in D_K$ we have
an estimate like \eqref{hast}, where 
$K'$ is a compact neighborhood of $K$.  However, in general
this estimate cannot be uniform in $x$ for $|x|\le 1$, since otherwise 
$\phi(a)$ would have a bounded extension to $X$. 

To see how this lack of uniformity may appear, assume that
 $a=f(b)$ for some hyperoperator $b$, where $f$ takes values
in $K\subset\subset \R^m$. Then  $D_{K,a}=D_b=D_a$ because if $x\in D_b$ 
then $\chi(b)x=x$ and since  $b$ has finite order on 
$\mbox{supp}\, \chi$ we get
\begin{multline*}
|\phi(a)x|=|(\phi\circ f \cdot \chi)(b)x|\leq
C\sum_{|\alpha|\leq N_{\chi}}\sup |\partial^{\alpha}(\phi\circ f \cdot \chi)||x|\\
\le C\sum_{|\alpha+\beta+\gamma|\leq N_{\chi}}\sup |\partial^{\beta}f||\partial^{\gamma}\chi|
\sup_{K}|\partial^{\alpha}\phi||x|.
\end{multline*}
However, 
$N_{\chi}$ and $\sup |\partial^{\beta}f||\partial^{\gamma}\chi|$ may blow up as 
$\chi \rightarrow 1$.
\end{remark}

\section{Spectrum of a hyperoperator}\label{spect}

We first recall

\begin{prop}
Suppose that $a=(a_1,\ldots, a_n)$ is a tuple of  bounded commuting operators with real
spectra and 
resolvents with temperate growths, and $A$ is the corresponding hyperoperator on $\R^n$.
Then $\supp A$ is equal to the (Taylor) spectrum of $a$.
\end{prop}

For a proof, see \cite{AS}. In view of this result the following definition
is natural.

\begin{df}
For   $A\in H_{\D(\R^n)}(X)$, the  spectrum $\sigma(A)$ is 
the support of $A$ as a distribution.
\end{df}

When $A$ is identified with the who
$a$ we often write $\sigma(a)$ instead of $\sigma(A)$.
Notice that $A(\phi)$ only depends on the values of $\phi$  in a small
neighborhood  of $\sigma(A)$. If the spectrum of $A$ is compact, then clearly
$A$ has a continuous extension to a multiplicative
mapping  $\E(\R^n)\to \L(X)$.
For such an $A$ and $f\in\E(\R^n)$, we have that
$f(a)x=\lim A(f\chi_N)x=A(f)x$ for $x\in D$,  and thus 
the closure of $f(a)$ is equal to the bounded operator $A(f)$. 
Applying to  the identity mapping $\xi\mapsto \xi$ on $\R^n$ 
we get

\begin{prop}
Suppose that $A\in H_{\D(\R^n)}(X)$ and $\sigma(A)$ is compact in $\R^n$. Then
the closure $\bar a$ of $a$ is bounded,  and $[a]=A$. 
Moreover,  $\sigma(A)$ coincides with the
Taylor spectrum of $\bar a$.
\end{prop}

If $f\in\E(\R^n)$   has its support 
in the complement of $\sigma(A)$, then $f(a)x=0$ for all $x\in D$,
so the closure of $f(a)$ is $0_X$.

\begin{df}
For a who $b=(b,D)$ we introduce the weak spectrum $\sigma_w(b)$ defined
as the intersection of all closed sets $F$ such that 
$\phi(b)x=0$ for all $x\in D$ and $\phi$ with support in $\R^n\setminus F$.
\end{df}

Thus a point $p$ is outside $\sigma_w(b)$ if and only if for all $\phi$ with support 
sufficiently close to $p$ we have $\phi (b)x=0$ for all $x\in D$.
It follows that if $b$ happens to be a hyperoperator then
$\sigma_w(b)=\sigma(b)$.
In particular, if $bx=0$ for all $x\in D$, then 
$\sigma_w(b)=\sigma(0)=\{0\}$.

\begin{prop}
Let $b=(b,D)$ be a who and let $f\in\E(\R^n,\R^m)$. Then
$
\sigma_w(f(b))=f(\sigma_w(b)).
$
\end{prop}

\begin{proof}
If    $h\in\D(\R^m)$ has its   support outside $f(\sigma(b))$, then
$h\circ f$ vanishes in a neighborhood of $\sigma(b)$ so
$(h\circ f)(b)x=0$ for $x\in D$, i.e,  by definition,
$h(f(b))x=0$.  This means that
$
\sigma_w(f(b))\subset f(\sigma_w(b)).
$
For the converse inclusion, take any point $p$ outside $\sigma_w(f(b))$
and let $h$ be a function identically equal to $1$ in a neighborhood of $p$ and with
support outside $\sigma_w(f(b))$. Then if $y\in f^{-1}(p)$ we have $h\circ f$ 
identically equal to $1$ in a neighborhood of $y$. Hence for any $\phi$ with support
in this neighborhood $\phi\cdot (h\circ f)=\phi$. Since $h$ has support outside 
$\sigma_w(f(b))$ we have $h\circ f (b)x =h(f(b))x=0$ for $x\in D$ and so
$\phi(b)x=\phi\cdot (h\circ f)(b)x=\phi(b)h\circ f(b)x=0$ for $x\in D$. Thus
$f^{-1}(p)\cap \sigma_w(b)=\emptyset$, i.e. $f(\sigma_w(b)) \subset \sigma_w(f(b))$.
\end{proof}

Noting that $\sigma_w(b)=\sigma(b)$ when
$b$ is a (strong) hyperoperator we immediately get

\begin{cor}
If $A\in H_{\D(\R^n)}(X)$ and $f\in\D(\R^n,\R^m)$, or 
$f\in\E(\R^n,\R^m)$ is proper, then
$
\sigma(f(a))=f(\sigma(a)).
$
\end{cor}

Since $\sigma_w(0_{D_A})=\sigma(0_X)=\{0\}$ we have

\begin{cor}\label{spekkorr2}
If $A\in H_{\D(\R^n)}(X)$ and  $f\in\E(\R^n,\R^m)$ and
$f(a)x=0$ for all $x\in D_A$, then
$
\sigma(a)\subset f^{-1}(0).
$
\end{cor}


It is not true in general that $f(\sigma(a))$ bounded implies that
$f(a)$ is bounded (if $f$ is neither proper nor compactly supported). 
For instance, take $f(\xi)=\sin \xi^m$ and $a\sim\xi$
on $X=H^1(\R)$. Then $|f|\le 1$ on $\sigma(a)$ but $f(a)$, i.e.,
multiplication with $\sin\xi^m$ is not bounded on $X$.
However we have

\begin{lma}
If $a$ is a hyperoperator, $f\in\E(\R^n)$,  and $b=f(a)$ is bounded, then
$f(\sigma(a))\subset \sigma(f(a))$.
\end{lma}

\begin{proof} We know that $p\circ f(a)=p(b)$ for all polynomials.
Let $\phi\in\D(\R^n)$ have support outside $\sigma(b)$ and take
$p_j$ such that $p_j\to \phi$ in $\E(\R^n)$. Then $p_j\to 0$ uniformly
in a neighborhood of $\sigma(b)$; we may even assume that this holds in a complex
neighborhood; thus we can conclude that $p_j(b)\to 0$ (even though we do not
know whether $b$ admits a smooth functional calculus or not!).

Moreover, $p_j\circ f\to \phi\circ f$ in $\E(\R^n)$ so 
$p_j\circ f(a)x\to \phi\circ f(a)x$ for $x\in D$. Since
$p_j\circ f(a)= p_j(b)\to 0$ we conclude that
$\phi\circ f(a)=0$. From Corollary \ref{spekkorr2} we get
$f(\sigma(a))\subset\{\phi=0\}$ and we conclude that
$f(\sigma(a))\subset \sigma(b)$.
\end{proof}

\begin{prop}\label{propp}
Assume that $a=(a,D)$ is a who   and that the closure
of $r_z(a)$ is bounded for each
$r_z(\xi)$, $z\in\C\setminus\R$. Then
$\bar a$ has real spectrum in the usual sense.
\end{prop}

\begin{proof}
We first prove that the closure $b$ of $r(a)=r_i(a)$ is the inverse of  $\bar a+i$.
We know that
$(a+i)b x=x=b(a+i)x$ for $x\in D.$
Suppose that $x\in \Dom(\bar a+i)=\Dom (\bar a)$. Then there are $x_j\in D$ such that
$x_j\to x$ and $(a+i)x_j\to (\bar{a}+i)x$. Since $b$ is bounded we have
$$
x\ot x_j=b(a+i)x_j\to b(\bar a+i)x
$$
so $b(\bar a+i)x=x$ for $x\in \Dom (\bar a+i)$. Moreover, if $x$ is arbitrary
and $x_j\in D$ and $x_j\to x$, then $b x_j\to bx$ and
$(a+i)bx_j=x_j\to x$ so by definition $bx$ is in the domain
of $\bar a+i$ and $(\bar a+i)bx=x$.
\end{proof}



\section{Representation by  pseudoresolvents}

We first consider the case $n=1$.
If $a$ is an HS operator, then we have the representation
\eqref{rep} of $\phi(a)$.
For a general $a\in H_{\D(\R)}(X)$ such a representation cannot hold
simply because  the resolvent is not defined.
We will discuss various ways to obtain formulas that will replace
\eqref{rep}.
The simplest way is  to use cut-off functions $\chi$ and define
$$
\omega_{\zeta-a}^\chi=\chi(\xi)\omega_{\zeta-\xi}|_{\xi=a},
\quad \omega_{\zeta-\xi}=d\zeta/(\zeta-\xi)2\pi i.
$$

\begin{prop}\label{intframst}
Suppose that $a\in H_{\D(\R)}(X)$. Then $\omega^\chi_{\zeta-a}$ is holomorphic for
$|\Im \zeta|>0$ and 
\begin{equation}\label{blik}
\| \omega^\chi_{\zeta-a}\|=\O(|\Im\zeta|^{-M})
\end{equation}
for some $M$.
If $\phi\in\D(\R)$ and $\supp\phi\subset\subset\{\chi=1\}$, then
\begin{equation}\label{eq3}
\phi(a)=\int\omega^\chi_{\zeta-a}\w\dbar\tilde\phi(\zeta).
\end{equation}
\end{prop}

\begin{proof}
By Lemma \ref{hololemma} $\omega_{z-a}^\chi$ is strongly holomorphic in 
$\C\setminus \R$.
Since $A$ has finite order on $K\supset\supset\supp\chi$, $A(\psi)$ only depends 
on a finite number of derivatives of $\psi$ if $\mbox{supp}\, \psi\subset K$  
and so we get \eqref{blik}.  If 
$$
\phi_{\epsilon}(\xi)=\frac{1}{2\pi i}
\int_{|\Im\zeta|>\epsilon}\frac{\chi(\xi)d\zeta}{\zeta-\xi}\w\dbar\tilde\phi(\zeta),
$$
it is readily checked, for instance by approximating by Riemann sums, that  
\begin{equation}\label{blik2}
\phi_\epsilon(a)=\int_{|\Im\zeta|>\epsilon}\omega^\chi_{\zeta-a}\w\dbar\tilde\phi(\zeta).
\end{equation}
Moreover, $\phi_{\epsilon}\to\phi\chi=\phi$ in $\D(\R)$, and hence
$\phi_{\epsilon}(a)\to\phi(a)$. Because of \eqref{blik}
it follows that the right hand side of \eqref{eq3} is absolutely  convergent 
and equal to the limit of the right hand side of \eqref{blik2}.
\end{proof}

\begin{prop}\label{spekkar}
Each $\omega^\chi_{z-a}$ has a holomorphic continuation to the set 
$\C\setminus\sigma(a)$;
more precisely, $\C\setminus\sigma(a)$
is  precisely the set where all $\omega^\chi_{\zeta-a}$ are strongly holomorphic.
\end{prop}

\begin{proof}
The first statement is proved analogously to Lemma \ref{hololemma}. If 
$x\in \R\setminus \sigma(a)$,  let $\tilde{\chi}$ be a cut-off function that is equal to 
$\chi$ in a neighborhood of $\sigma(a)$ and zero in a neighborhood of $x$.
Then $\omega^\chi_{z-a}=A(\tilde{\chi}/(z-\xi))$ and imitating the proof of 
Lemma \ref{hololemma} we see that $\omega^\chi_{z-a}$ is strongly holomorphic close
to $x$. 
For the converse, assume $\phi$ has its support where $\omega^\chi_{\zeta-a}$ 
is holomorphic and $\chi$ identically $1$ in a neighborhood  of $\supp \phi$.
Then by Proposition~\ref{intframst},
$$
A(\phi)=\int\omega^\chi_{z-a}\w\dbar\tilde\phi(z)=
-\int\dbar(\tilde{\phi}(z)\omega^\chi_{z-a})=0
$$
by Stokes'  theorem and thus we are done.
\end{proof}

The advantage with the usual representation  \eqref{rep} is of course that
a~priori we only have to compute $\phi(a)$ for $\phi(\xi)=1/(\zeta-\xi)$.
For the  general hyperoperator  we must insert various  functions  $\chi$ as well.
However, if we impose growth restrictions on $[a]$, one single formula  will do.
In Section~\ref{thyp} we will consider the case with polynomial
growth restrictions.

\smallskip
If $a$ is a hyperoperator or even just a who, then for each
$x\in D$,  the resolvent $\omega_{\zeta-a}x$ is holomorphic outside
the compact set $\sigma_x(a)\subset\R$, and  from  \eqref{hast} we have that
$|\omega_{\zeta-a}x|\le C|\Im\zeta|^{-M}$. With a similar argument as above we therefore have
the representation
$$
\phi(a)x=\int\omega_{\zeta-a}x\w\dbar\tilde\phi(\zeta),\quad  \phi\in\E(\R).
$$

Recall that $\G(\mathbb{R})$ is the algebra of functions on 
$\widehat\R$ that are holomorphic in a complex
neighborhood of $\infty$.
Convergence in $\G(\R)$ of a sequence $f_j$ means that  $f_j$ converges in
$\E(\widehat\R)$ and moreover, that all $f_j$  are holomorphic in a fixed complex 
neighborhood of  $\infty$ and converge uniformly on compacts in this neighborhood.
\begin{thm}\label{gsats}
A hyperoperator  $A\in H_{\D(\R)}(X)$ 
corresponds to an HS operator  if and only if
$A\colon\D(\R)\to\L(X)$ has a multiplicative continuous extension to a mapping
$\G(\R)\to\L(X)$. 
\end{thm}


This result was more or less proved in  \cite{AS}; one part  is contained in the proof of 
 Proposition $7.2$ in \cite{AS} and the other part is stated in Proposition $11.4$ in 
the same paper, but for the reader's convenience we supply a proof here.

\begin{proof}
First we notice that such an extension of $A$ must be unique if it exists at all.
In fact,  for any $\psi \in 
\G(\mathbb{R})$ and $x\in D_A$ we have
$
A(\psi \chi)x=\psi(a)x
$
if $\chi$ is chosen so that $A(\chi)x=x$. On the other hand if $\hat{A}$ is a 
multiplicative extension of $A$ we get
$
\hat{A}(\psi)x=\hat{A}(\psi)A(\chi)x=A(\psi \chi)x
$
Hence $\hat{A}(\psi)$ coincides with $\psi(a)$ on $D_A$ and since $D_A$ is dense and 
$\hat{A}(\psi)$ is bounded this uniquely determines $\hat{A}(\psi)$. 
Here  $a$ denotes the who  associated to $A$.

For the ``only if''-part we first assume that (the closure of) $a$ 
is an HS operator, cf., Proposition~\ref{sallad3}.
Then the action of $A$ is given by \eqref{rep} and we
want to extend this formula
to any function $f$ in $\G(\R)$. Let $F$ be the holomorphic extension
to a complex neighborhood $O$ of $\infty$, and let $\chi$ be a cut-off function in 
$\R$ that is equal to $1$ in a neighborhood of $K=\R\setminus (\R\cap O)$.
One can find an almost holomorphic extension $\tilde\chi$ which is
$0$ in a complex neighborhood of $\infty$ and $1$ in a complex neighborhood 
of $K$. Then 
$
\tilde f=(1-\tilde\chi)F+\widetilde{\chi f}
$
is an almost holomorphic extension of $f$ to a complex neighborhood  of $\widehat{\R}$ 
in $\widehat\C$ 
which is holomorphic in a neighborhood of $\infty$. 
Let $\psi$ be a function identically equal
to $1$ in a neighborhood of $\widehat{\R}$ in $\widehat{\C}$ and with support in a slightly 
larger neighborhood avoiding the point $i$.  Then
\begin{equation}\label{eqmojs}
\frac{1}{2\pi i}\int\frac{a-i}{\zeta-i}\frac{d\zeta}{\zeta-a}\w\dbar(\tilde{f}\psi
(\zeta))
\end{equation}
provides the desired extension. In fact, if $f$ has compact support then $\tilde{f}
\psi$ is an almost holomorphic extension of $f$ with compact support avoiding $i$.
It follows by Stokes' theorem  that formula \eqref{eqmojs} yields the same operator
as \eqref{rep}. Moreover \eqref{eqmojs} is continuous and multiplicative on 
$\mathcal{G}(\R)$. This is perhaps most easily seen by pulling back to the unit circle
$\mathbb{T}$. The Cayley transform  $b=C(a)$,
cf., \eqref{cayley}, is  a bounded operator with spectrum contained in
the unit circle $\T$, and
$$
\frac{1}{2\pi i}\int\frac{a-i}{\zeta-i}\frac{d\zeta}{\zeta-a}\w\dbar(
\tilde{f}\psi(\zeta))=
\frac{1}{2\pi i}\int\frac{dw}{w-b}\wedge \dbar (\tilde{f}\psi(C^{-1}(w))).
$$
The right hand side is a continuous extension of the holomorphic functional calculus 
for $b$ to the space of smooth functions on $\mathbb{T}$ which are analytic
in a neighborhood of $1$ since $\|(w-b)^{-1}\|$ has tempered growth in 
$\mathbb{T}\setminus \{1\}$. Since the analytic functions are dense in this space,
the multiplicativity follows automatically.

Conversely, assuming that  $A$ is a hyperoperator that admits  an extension to $\G(\R)$,
we want to prove that $\bar a$ is an HS operator.
Since $A$ now operates on all $r_z(\xi)=1/(z-\xi)$ it follows 
from Proposition \ref{propp} that $\bar{a}$ has spectrum in $\R$ in the usual sense.
Clearly then $r_z(a)dz/2\pi i$ is the resolvent of $\bar{a}$.
Given a compact $K\subset\R$ take $\chi$ and $\tilde\chi$ as above.
As $A$ has finite order $m$ on $K'=\mbox{supp}\,\chi$ it follows that 
\begin{equation}\label{brunn}
\Big\|  \frac{\chi(a)}{z-a}\Big\|\leq C_{K} |\Im z|^{-(m+1)}
\end{equation}
for any $z\in \C\setminus \R$. 
For $z$ in a small neighborhood of $K$, the functions 
$$
g_z(\xi)=\frac{\tilde\chi(z)-\chi(\xi)}{z-\xi}.
$$
are uniformly bounded in  $\G(\R)$, and by \eqref{brunn} so are
$\| g_z(a)\|=\| A(g_z)\|$. Thus  \eqref{tempupp} 
follows by the triangle inequality.
\end{proof}

\begin{remark}\label{rem10}
Let $a\in H_{\D(\C)}(X)$. Then we can define
 $\omega^\chi_{\zeta-a}$ as a $\L(X)$-valued distribution ($(1,0)$-current) in $\C$ by
\begin{equation}\label{abstrakt}
\omega^\chi_{\zeta-a}.\psi d \bar\zeta=
\frac{1}{2\pi i}\int_\zeta\frac{\chi(\xi)\psi(\zeta)d\zeta\w d\bar\zeta}{\zeta-\xi}\Big|_{\xi=a},\quad \psi\in\D(\C).
\end{equation}
If we apply to $\dbar\psi$ we get
$
\omega^\chi_{\zeta-a}.\dbar\psi=\chi(\xi)\psi(\xi)|_{\xi=a}=\chi(a)\psi(a).
$
Thus 
$\dbar \omega^\chi_{\zeta-a}=\chi[a].$
If in fact $a\in H_{\D(\R)}(X)$ and
we  choose $\psi=\tilde\phi$ in \eqref{abstrakt}, then we can move  $a$ inside  the integral
and thus  get back  \eqref{eq3}.
However, in general it is not possible to  put $a$ inside the  integral.
\end{remark}

If we want an absolutely convergent integral representation for
$\phi(a)$ when $a\in H_{\D(\R^n)}(X)$ we can use the  Bochner-Martinelli form
$$
\omega_\xi=b(\xi)\wedge (\dbar b(\xi))^{n-1},\quad b(\xi)=
\frac{\sum \bar\xi_jd\xi_j}{2\pi i|\xi|^2},
$$
and define  
$
\omega^{\chi}_{\zeta-a}=\chi(\xi)\omega_{\zeta-\xi}|_{\xi=a}.
$
Then $\omega_{\zeta-a}$ is $\dbar$-closed in $\C^n\setminus\R^n$ and
the analogue of Proposition \ref{intframst} holds.
Proposition \ref{spekkar} also has a generalization to the $\R^n$ case; $\C^n 
\setminus \sigma(a)$ is precisely the set where $\omega^{\chi}_{z-a}$ is strongly
$\dbar$-closed.
If we consider a hyperoperator $a\in H_{\D(\R^{2n})}(X)$ as an element
in  $a\in H_{\D(\C^{n})}(X)$, the analog of Remark~\ref{rem10} also holds.
\smallskip

Tensor products of hyperoperators  can also be defined by integral formulas.
Assume that $A_1,\ldots,A_m$ are in $H_{\D(\R^{n_j})}$ but not necessarily
commuting. Then we can form the tensor product $A=A_1\otimes\cdots \otimes A_m$,
and obtain a linear continuous, though not multiplicative,
operator $\D(\R^n)\to\L(X)$, where $n=n_1+\cdots +n_m$.
For $\phi\in\D(\R^n)$ we can find an almost holomorphic extension $\tilde\phi$ such that
\eqref{asspecial} holds. In \cite{AS} this is only proved when all $n_j=1$ but the 
general case follows  along the same lines.
Then
\begin{equation}\label{defa}
(A_1\otimes\cdots\otimes A_m)(\phi)
=\int\omega^{\chi_1}_{\zeta_1-a_1}\w\ldots
\w\omega^{\chi_m}_{\zeta_m-a_m}\w\dbar_{\zeta_m}\cdots\dbar_{\zeta_1}\tilde\phi(\zeta),
\end{equation}
if the support of $\phi$ is contained in the set where $\chi_1\otimes\cdots\otimes\chi_m=1$. 
To see this, first notice that the integral makes sense in view of the assumption
\eqref{asspecial} and the estimates \eqref{blik} of $\omega^{\chi_j}_{\zeta_j-a_j}$. Since 
\eqref{defa} clearly holds for 
$\phi$ of the form $\phi=\phi_1\otimes\cdots\otimes\phi_m$,  the general case follows by
continuity. 
One can also prove directly that \eqref{defa} is independent of the choice of special
almost analytic extension $\tilde\phi$ along the lines in \cite{AS}, and then use this
as the definition of the tensor product.

\begin{remark}
We can also generalize Theorem \ref{gsats} to several variables, and
we  illustrate by considering  a hyperoperator 
$A\in H_{\D(\R^2)}(X)$. 
First we define $\mathcal{G}(\R^2)$ as the union (direct limit) of the spaces
$\mathcal{G}_U(\R^2)$, $U$ a complex neighborhood of $\infty$ in $\widehat{\C}$, defined as 
all smooth functions $f$ on  $\widehat{\R}\times \widehat{\R}$ which are holomorphic 
on $U\times U$ and such that $x\mapsto f(x,y)$ is holomorphic  in $U$ for any  $y$
and $y\mapsto f(x,y)$ is holomorphic in $U$ for any $x$.  A sequence $f_j$ in 
$\mathcal{G}(\R^2)$ converges if all $f_j$ are in some fixed $\mathcal{G}_U(\R^2)$
and converges in $\E((\widehat{\R}\cup U)\times (\widehat{\R}\cup U))$.
The analog of Theorem 
\ref{gsats} is:   $A$ has a continuous extension to $\mathcal{G}(\R^2)$ if and only if 
the closures  of $a_j=A(\pi_j)$, $j=1,2$, are of HS type and commute strongly,
i.e., their resolvents commute.
Notice however that this  condition highly depends on the choice of
coordinates on $\R^2$, whereas the notion of general  hyperoperator
is coordinate invariant.
\end{remark}

\section{Temperate hyperoperators}\label{thyp}
We say that $A\in H_{\D(\R^n)}(X)$ is temperate, $A\in H_{\S(\R^n)}(X)$, if
it extends to a (necessarily multiplicative)
mapping $\S(\R^n)\to \L(X)$.

\smallskip
Since $\D(\R^n)$ is dense in $\S(\R^n)$ it follows that a continuous
multiplicative map $\S(\R^n)\to \L(X)$ satisfies $(i)$ and $(ii)$ in 
Definition~\ref{hodef} 
if and only if it holds with $\D(\R^n)$ replaced by $\S(\R^n)$
(but the corresponding dense domain may be larger).

\smallskip

For standard functional analysis reasons it follows that
for any temperate $A$ there is an integer $M$ such that
\begin{equation}\label{eq15}
|A(\phi)|\le C \sum_{|\alpha|,|\beta|\le M}\sup_{\R^n} |\xi^\beta\partial^\alpha \phi|,
\end{equation}
which in particular means that $A(\phi)$ is defined for $\phi$
such that its derivatives up to order $M$ as least have decay like
$1/|\xi|^M$.

\begin{ex}
Let $X$ be the set of functions $\phi(\xi)$  on $\R$ with norm
$\|\phi\|=\sum_\ell \|\phi\|_{C^\ell( K_{\ell+1}\setminus int K_{\ell-1})}$.
Then multiplication with $\xi(2+\sin \xi^3)$ is a 
hyperoperator that is not temperate.
\end{ex}

The multiplication hyperoperator $f(\xi)=\xi(2+\sin\xi^3)$
on $H^1(\R)$ from Example \ref{ex100} is a tempered hyperoperator, which has no ordinary 
resolvent. Notice, though, that
$$
\frac{i+\zeta}{i+f(\xi)} \frac{1}{\zeta-f(\xi)}
$$
is bounded for all $\zeta\in\C\setminus\R$. 
More generally,
if    $A\in H_{\S(\R)}(X)$ and $m$ is a large enough integer we can define, 
in view of \eqref{eq15},
$$
\omega^m_{\zeta-a}=\left(\frac{i+\zeta}{i+\xi}\right)^m\omega_{\zeta-\xi}\big|_{\xi=a},
$$
for $\zeta\in\C\setminus\R$. 
If  $A\in H_{\S(\R^n)}(X)$ we can take instead 
$$
\omega^m_{\zeta-a}=\left(\frac{1+\zeta\cdot\xi}{1+|\xi|^2}\right)^m\omega_{\zeta-\xi}\big|_{\xi=a}.
$$
for $\zeta\in\C^n\setminus\R^n$.

\begin{prop}
The  form $\omega^m_{\zeta-a}$ is $\dbar$-closed in  $\C^n\setminus\R^n$ and admits a
$\dbar$-closed  extension
to $\C\setminus\sigma(a)$. 
Moreover,
if $\phi\in\S(\R^n)$ and $\tilde\phi$ is an appropriate almost holomorphic extension, then
\begin{equation}\label{eq16}
A(\phi)=\int \omega^m_{\zeta-a}\w\dbar\tilde\phi.
\end{equation}
\end{prop}

This means,  cf., Remark~\ref{rem10}, that $\dbar\omega^m_{\zeta-a} =[a]$.
Moreover, if $\omega^m_{\zeta-a}$ has a $\dbar$-closed extension to $\C^n\setminus F$, then
$\sigma(a)\subset F$.

\begin{proof}[Sketch of proof]
First notice that  
$$
\sup_{\xi\in\R^n}|\xi^\beta\partial^\alpha_\xi \omega^m_{\zeta-\xi}|\le 
C\frac{(1+|\zeta|)^m}{|\Im\zeta|^{2n-1+|\alpha|}}
$$
if just $|\beta|<m$. If $A$ satisfies \eqref{eq15}, therefore  $\omega^m_{\zeta-a}$ is 
well-defined if $m\ge M$, and
\begin{equation}\label{eq17}
\|\omega^m_{\zeta-a}\|\le C\frac{(1+|\zeta|)^m}{|\Im\zeta|^{2n-1+|\alpha|}}.
\end{equation}
Given $\phi\in\S(\R^n)$ we let
$$
\tilde \phi(\zeta)=\int_t e^{it\cdot\zeta}\hat \phi(t)\chi\big(\sqrt{1+|t|^2}|\Im\zeta|\big),
$$
where $\chi(s)$ smooth, supported in the unit ball in $\R^n$ and
 identically $1$ in a neighborhood of the origin. One easily checks that
$\tilde \phi(\zeta)$ is smooth, and equal to $\phi$ on $\R^n$, and that moreover,
\begin{equation}\label{astrid}
\dbar \tilde{\phi}(\zeta)=
\O_{M_1,M_2}(|\Im\zeta|^{M_1}(1+|\zeta|)^{-M_2}), \quad M_1,M_2>0.
\end{equation}
In view of \eqref{eq17}, therefore,  the integral in \eqref{eq16} is well-defined.
Moreover, from \eqref{astrid} it is easily seen  that 
$
\int\omega_{\zeta-\xi}\w\dbar\tilde\phi(\zeta)=\phi(\xi),
$
and replacing $\tilde\phi(\zeta)$ by 
$$
\tilde\phi(\zeta)\Big(\frac{1+\zeta \cdot \xi}{1+|\xi|^2}\Big)^m
$$
which satisfies a similar estimate, we get that
$$
\int\omega^m_{\zeta-\xi}\w\dbar\tilde\phi(\zeta)=\phi(\xi).
$$
One then proves \eqref{eq16}  along the same lines as Proposition~\ref{intframst}.
\end{proof}

For tempered hyperoperators  the theory for tempered distributions is at our disposal.
We will use this to prove a new form of Stone's  theorem.
We first recall a simple known variant. 

\begin{ex}
If $v\in C^1(\R^n,\L(X))$ and
\begin{equation}\label{krut}
v(t+s)=v(t)v(s), \quad v(0)=e_X,
\end{equation}
then $v(t)=e^{ia\cdot t}$, for the commuting tuple  $a_k=(\partial v/\partial t_k)(0)/i$ 
in $\L(X)$. If in addition $|v(t)|=\O(|t|^m)$, when $|t|\to\infty$, then 
$\sigma(a)\subset\R^n$.
If we only assume that $v(t)$ is continuous and satisfies \eqref{krut}, then the conclusion is not true.
(For instance,  if  $n=1$ and $a$ is multiplication with $\xi$ on $L^2(\R_\xi)$, then
$v(t)=e^{iat}$ is multiplication by $e^{i\xi t}$ and thus
a bounded operator, but  $v'(0)$ is not bounded.)  
However,  $v$ is generated by a hyperoperator  $A\in H_{\S(\R^n)}(X)$, i.e.,
$v(t)=\exp ia\cdot t$.

In fact, 
assume that  $v(t)$ is continuous in the weak sense that $v(t)x$ is continuous for each
$x\in X$. It  then  follows from the  Banach-Steinhaus  theorem that 
$\|v(t)\|$ is uniformly bounded on compact sets.
Therefore, $v.\phi=\int_t v(t)\phi(t)dt$
is a bounded operator for each $\phi\in\S(\R^n)$.
Moreover, the condition $v(0)=e_X$ implies that
$\cap \Ker v(\phi)=\{0\}$ and $\cup \Im v(\phi)$ is dense. 
In fact, let $\phi_j\to\delta_0$.
Then $v(\phi_j)x\to x$ since $\varphi\mapsto v(\varphi)x$ is continuous and we easily see
that $\cap \Ker v(\phi)=\{0\}$ and that $\cup \Im v(\phi)$ is dense. 
The existence of the generator $A$ now follows from
Proposition~\ref{stone} below.
\end{ex}

Let  $A$ be a tempered hyperoperator and let 
$
D=\cup_{\phi\in\S(\R^n)}\Im A(\phi).
$
If  $f\in\E(\R^n)$  is a multiplier
on $\S(\R^n)$, i.e., $f\S(\R^n)\subset\S(\R^n)$, we can define
$f(a)x$ for $x\in D$ as $A(f\phi)y$ if $x=A(\phi)y$. To see that this is
well-defined, assume that also $x=A(\phi')y'$. By the multiplicativity,
we then have  that $A(\chi_Nf\phi)y=A(\chi_Nf\phi')y'$ since
$\chi_Nf$ is in $\S$. When $N\to\infty$, $\chi_Nf\phi\to f\phi$ in $\S$,
and hence $A(f\phi)y=A(f\phi')y'$.
It is readily checked that $f(a)$ maps $D\to D$ and that $(fg)(a)x=f(a)g(a)x$.

Observe that $f(\xi)=\exp(i\xi\cdot t)$ is a multiplier on $\S$, so
$\exp(ia\cdot t)x$ is defined for all $x\in D$. Moreover,
 $x=A(\phi)y$ so $\exp(ia\cdot t)x=A(\phi(\xi)\exp(i\xi\cdot t))y$, 
and therefore     \eqref{eq15} implies that 
\begin{equation}\label{baster}
|e^{ia\cdot t}x|\le C_x|t|^M.
\end{equation}
We claim that 
\begin{equation}\label{skoda}
A(\hat\psi)x=\int_t\psi(t)e^{-ia\cdot t}x dt, \quad \psi\in\S, x\in D.
\end{equation} 
In fact, the integral is convergent in view of \eqref{baster} and
it is easy to see that it is equal to $A(\hat\psi)x$ since
$
\int_{|t|<R}\psi(t)e^{-i\xi\cdot t}dt\to\hat\psi(\xi)
$
in $\S$. In particular, the integral in \eqref{skoda} has a continuous
extension to $X$. Since $A$ is in $\S'$ it has a Fourier transform $\hat A$, 
defined by $\hat A(\psi)=A(\hat\psi)$,  and thus we have the suggestive formula
$\hat A(t)=\exp(-ia\cdot t)$.
If we let $v(t)=\exp(-ia\cdot t)=\hat A(t)$ then clearly
$v(t+s)x=v(t)v(s)x$ for $x\in D.$
Moreover, clearly $v$ defined by $v.\psi=\hat\psi(a)$ satisfies
\begin{equation}\label{full1}
\int_s\int_t v(t+s)\phi(t)\psi(s)=\int_tv(t)\phi(t)\int_sv(s)\psi(s),\quad \phi,\psi\in\S
\end{equation}
and 
\begin{equation}\label{full2}
\cap \Ker v(\phi)=\{0\},\quad \cup \Im v(\phi)=D\ {\rm dense}.
\end{equation}
We have the following variant of Stone's theorem.

\begin{prop}\label{stone}
Assume that $v\colon \S(\R^n)\to\L(X)$ is linear, and continuous in the sense that
for fixed $x\in X$, $v.\phi_j x\to 0$ whenever $\phi_j\to 0$ in $\S(\R^n)$.
Moreover, assume that $v(t)$ is  group of operators  in the sense of \eqref{full1}
and $v(0)=e_X$ in the sense of \eqref{full2}.
Then $v$ is generated  by a hyperoperator $A\in H_{\S(\R^n)}(X)$ in the sense that
$v(t) x$ is smooth for $x\in D_A$ and $(\partial v/\partial t_k)(0)x=ia_k x$.
\end{prop}

\begin{proof}
Define $A(\phi)=v.\hat{\phi}$. By the Banach-Steinhaus theorem the pointwise continuity 
of $v$ implies strong continuity and so $A$ is a continuous map $\mathcal{S}(\R^n)\rightarrow
\mathcal{L}(X)$. Moreover, the weak multiplicativity of $v$ implies that $v(\phi*\psi)
=v(\phi)v(\psi)$ and hence
\[
A(\phi \psi)=v(\widehat{\phi \psi})=v(\hat{\phi}*\hat{\psi})=v(\hat{\phi})v(\hat{\psi})=
A(\phi)A(\psi).
\]
Since the Fourier transform is an isomorphism of $\mathcal{S}(\R^n)$ we get that
$\cap \Ker A(\phi)=\{0\}$ and $\cup \Im A(\phi)=D$ is dense. Thus $A$ is a tempered 
hyperoperator. 
For $x\in D$ we can define $u(t)x=e^{iat}x$ and since $A$ satisfies an estimate like
\eqref{eq15} it is easy to see that $|u(t)x|\leq C|t|^M$ and so $u(t)x$ defines an 
element in $\mathcal{S}'(\R^n,X)$. We also see that $t\mapsto u(t)x$ is in $C^1$ (even
in $C^{\infty}$) and $u'(t)x=iax$. In fact, if $\phi \in \mathcal{S}$ then 
$\phi(\xi)(e^{i\xi t}-1)/t\rightarrow i\xi \phi(\xi)$ in $\mathcal{S}$
as $t\rightarrow 0$, and hence if $x=A(\phi)y$ we get
\[
\frac{e^{iat}-e_X}{t}x=\frac{\phi(a)e^{iat}-\phi(a)}{t}y\rightarrow
ia\phi(a)y=iax.
\]
We finally check that $u(t)x=v(t)x$ as tempered distributions. If,   as before,  $x=A(\phi)y$,
then for any $\psi \in \mathcal{S}$ we have
\begin{eqnarray*}
\int_t\psi(t)u(t)x &=& \int_t \psi(t)A_{\xi}(\phi(\xi)e^{i\xi t})y=
A_{\xi}(\phi(\xi)\int_t \psi(t)e^{i\xi t})y \\
&=&
A_{\xi}(\phi(\xi)\hat{\psi}(\xi))y
=
A(\hat{\psi})x=
\int_t \psi(t)v(t)x.
\end{eqnarray*}
\end{proof}

\section{Operators with ultradifferentiable functional calculus}\label{ultra}

Let $h(t)=H(|t|)$ where $H(0)=0$ and $H$ increasing and concave
on $[0,\infty)$.
Then $h$ is subadditive. We also assume that
$\lim_{|t|\to\infty}h(t)/|t|=0$
and that
\begin{equation}\label{limsup}
\limsup_{|t|\to \infty}\frac{\log(1+|t|)}{h(t)}=0.
\end{equation}
Let
$\A_h$ be the space of  tempered  distributions  $f$ on $\R^n$ such that
$\hat f$ is a measure and
\begin{equation}\label{hatest}
\|f\|_{\A_h}=\int_t|\hat f(t)|e^{h(t)}dt < \infty.
\end{equation}
Because of   \eqref{limsup}, $\A_h$ is contained in $C^{\infty}(\R^n)$. 
Clearly  $\A_h$ is a Banach space of functions that is 
closed  under translations, and since $h$ is subadditive it follows,
see e.g., \cite{AB}, that
$\A_h$ actually is a Banach algebra under pointwise multiplication. 
These algebras were introduced by Beurling, \cite{Beu}. If
$h(t)=|t|^{\alpha}$, $0<\alpha<1$, then
$G_\alpha=\cup_{c>0}\A_{ch}$ is the classical Gevrey algebra, see \cite{H}.
We say that the class $\A_h$ is  non-quasianalytic if for each
compact set $E$ and open neighborhood $ U\supset E$  there is a
function $\chi\in\A_h$ with  support in
$U$ which is  identically $1$ in some neighborhood of $E$.
We  recall the following version of the
Denjoy-Carleman theorem. 

\begin{thm}\label{HB}
The class $\A_h$  is non-quasianalytic  if and only if
\begin{equation}\label{Hvillkor}
\int_1^{\infty}\frac{H(s)ds}{s^2}<\infty.
\end{equation}
\end{thm}

Assume now that $h(t)=H(|t|)$ satisfies the condition \eqref{Hvillkor}.
Let $B_{h}$ be the algebra of all functions on  $\R^n$ which are locally
in $\A_{ch}$ for some $c>1$, and let $B_{h,0}$ be the subalgebra
of functions with compact support.
There is an associated convex decreasing function
$G(s)=\sup_t(H(t)-ts)$ on $(0,\infty)$.
Let $H_c(s)=H(cs)$ and let $G_c$ be the corresponding decreasing function.

\begin{prop}\label{ek}
A function $\phi\in B_h$ if and only if it admits an almost holomorphic extension
$\tilde\phi$ such that for each compact $K\subset\R^n$, for some $c>1$ we have
$$
\sup_{\Re\zeta\in K}|\dbar\tilde\phi|e^{g_c(\Im\zeta)}<\infty.
$$
If $\phi$ has compact support and $U$ is a complex neighborhood of
$\supp\phi$ we can choose $\tilde\phi$ with support  in $U$.
\end{prop}

For a proof, see, e.g., \cite{AB}.
It follows that composition of functions in $B_h$ stays in $B_h$.
In a completely analogous way as before we can now define
a hyperoperator $A\in H_{B_{h,0}}(X)$ as a continuous multiplicative mapping
$B_{h,0}(\R^n)\to \L(X)$
such that
$\cup_{\phi\in B_{h,0}}\Im A(\phi)=D$ is dense and
$\cap_{\phi\in B_{h,0}}\Ker A(\phi)=\{0\}.$
Everything that is done in Sections~3,4, and 5 carry over directly to these
ultrahyperoperators;  for instance,
$D$ is the set of $x\in X$ such that $x=A(\chi)x$ for some
cut-off function $\chi$ in $B_h$. 
If $A\in H_{B_{h,0}}(X)$,  then 
$\|A(\phi)\|\le C_c\sup_{K'_c}|\phi|_{\A_{ch}}$
for each   $c>1$. If we define 
$\omega^\chi_{\zeta-z}=\chi(\xi)\omega_{\zeta-\xi}|_{\xi=a}$ it turns out that
$
\|\omega^\chi_{\zeta-z}\|\le C_c\exp{g_c(\Im\zeta)}
$
for each $c>1$. If $\supp\phi\subset\{\chi=1\}$
we thus have the representation 
$$
A(\phi)=\int\omega^\chi_{\zeta-a}\w\dbar\tilde\phi(\zeta).
$$

\section{Invariant subspaces and spectral decomposition}

Precisely as for a bounded operator (tuple of commuting bounded operators) that admits
a smooth functional calculus, for a hyperoperator $a$ there is a rich structure
of invariant subspaces as well as spectral decompositions.

\begin{prop}\label{nollrum}
Assume that $A\in H_{\D(\R^n)}(X)$, $f\in\E(\R^n,\R^m)$, and let
$$
X'=\{x\in D_A;\ f(a)x=0\}.
$$
Then $Y=\overline{X'}$ is an $a$-invariant subspace of $X$, and 
$a'=a|_Y$ is a hyperoperator. Moreover, $D_{a'}=X'$ and
\begin{equation}\label{ink}
{\rm int\, }\{f=0\}\cap\sigma(a)\subset\sigma(a')\subset\{f=0\}\cap\sigma(a).
\end{equation}
If $\{f=0\}$ contains some open subset of $\sigma(a)$, then 
$Y$ has nontrivial vectors.
\end{prop}

\begin{proof}
Since $f(a)$ and $\phi(a)$ commute, $X'$ and hence $Y$ are $a$-invariant.
If $\phi$ has compact support, then $\phi(a)$ is bounded,
and hence $\phi(a')$ extends to a bounded operator on $\overline{X'}$. Moreover, the
continuity with respect to $\phi$ is clear. Since $1(a')x=x$ for all $x\in X'$, the 
properties $(i)$ and $(ii)$ in Definition \ref{hodef} are  satisfied, 
so $a'$ is indeed a hyperoperator on the Banach space $Y$.

By definition, $X'\subset D_a$. If $x\in D_{a'}$, then $x=\chi(a')x$ 
for some $x\in Y$.
This means that $x=\chi(a)x$ and so $x\in D_a$, and moreover $f(a)x=0$. Thus $D_{a'}=X'$.

If $\phi(a)=0$ for all $\phi\in\D(\omega)$ then $\phi(a')=0$ for all such $\phi$, and hence
$\sigma(a')\subset\sigma(a)$. If $p$ is any point  outside $\{f=0\}$ 
then $f_j(p)\neq 0$ for
some $f_j$ ($f=(f_1,\ldots,f_m)$). We may assume that $f_j(p)=1$. 
If $\omega\ni p$ is small enough,
$|f_j-1|\le 1/2$ in $\omega$. For $\phi\in\D(\omega)$ we have that
$$
\phi(a)x=\phi(a)(1-f_j)^N(a)x=(\phi(1-f_j)^N)(a)x, \quad x\in X',
$$
and since $\phi(1-f_j)^N\to 0$ in $\D(\omega)$ when $N\to\infty$ we 
can conclude that $\phi(a)x=0$.
Thus $\omega$ is contained in the complement of $\sigma(a')$ and so 
 we have proved the second inclusion in \eqref{ink}
To see the first one, take $p\in{\rm int\, }\{f=0\}\cap\sigma(a)$ and a neighborhood $\omega$ such that 
$p\in\omega\subset\{f=0\}$. Since $\omega$ intersects $\sigma(a)$ there exists some $\phi\in\D(\omega)$
and $z\in X$ such that $x=\phi(a)z\neq 0$. However, then $x\in D$ and 
$f(a)x=(f\phi)(a)z=0$ since $f\phi=0$, so $x\in X'$.  Thus $\omega$ intersects $\sigma(a')$.
Since $\omega\ni p$ can be chosen arbitrarily small, we conclude that $p\in\sigma(a')$.
If $\sigma(a')$ is nonempty, then $Y$ is nontrivial, and so the last statement
follows from \eqref{ink}.
\end{proof}

If $p$ is an isolated point in $\sigma(a)$ and $f=0$ in a neighborhood of $p$, then
$X'$ is non-trivial. 
There are also non-trivial $a$-invariant subspaces as soon as $\sigma(a)$ contains more than
one point.
Notice that $a'$ is bounded if $\{ f=0\}\cap\sigma(a)$ is compact.

It is  easy to make spectral decompositions.
Let $A\in H_{\D(\R^n)}(X)$ be a hyperoperator and 
let $\{\Omega_j\}$ be a locally finite open cover of $\sigma(a)$.
Moreover, choose  $\phi_j\in\E(\R^n)$ such that
$\Omega_j\subset\{\phi_j=1\}$, and let
$$
X_j=\{x\in D_A;\ \phi_j(a)x=x\}.
$$
If $\Omega_j$ is bounded, we can choose $\phi_j$ in $\D(\R^n)$ and then
$X_j=\Ker (e_X-A(\phi))$ is a closed subspace of $D_A$.
Then $X_j$ are $a$-invariant subspaces, 
 $\sigma(a|_{X_j})\subset \Omega_j\cap\sigma(a)$,
and
\begin{equation}\label{uppdel}
\sum_1^\infty X_j=D_A.
\end{equation}
All these statements but the last one follows from  Proposition~\ref{nollrum}.
To see \eqref{uppdel}, choose a smooth partition of unity $\chi_j$ subordinate to
the cover $\{\Omega_j\}$. Then, since $\sum\chi_j=1$,
for each $x\in D_A$ we have 
$
x=\sum_1^M\chi_j(a)x
$
for some $M$. However, $(1-\phi_j)\chi_j=0$ so $\chi_j(a)x$ belongs to $X_j$.
Hence, \eqref{uppdel} follows.

In general the sum \eqref{uppdel} is not direct. However, if 
 $\sigma(a)$ is a disjoint union of closed sets $F_j$, we can find $\phi_j$
with disjoint supports such that $\{\phi_j=1\}$ contain a neighborhood of $F_j$.
If $x\in X_j\cap X_k$, then
$x=\phi_j(a)x=\phi_j(a)\phi_k(a)x=(\phi_j\phi_k)(a)x=0$, and hence we get
$$
D_A=\oplus_1^\infty X_j.
$$

\begin{ex}
Let $A\in H_{\D(\R^n)}(X)$ and let   $f\in\E(\R^n,\R^m)$ be a mapping such that
 $f(a)=0$. 
From Corollary \ref{spekkorr2} (or \eqref{ink}) we know   that
$\sigma(a)\subset\{f=0\}$.   Let us also assume that the zero set $\{ f=0\}=\{\alpha^j\}$ 
is discrete. Then we have the decomposition
$
D_A=\oplus_1^\infty X_j
$
where $\sigma(a|_{X_j})=\{\alpha^j\}$.
For each $j$, let $g_1^j,\ldots, g^j_{\ell_j}$ be functions in the local ideal generated
by $f$ at $\alpha^j$, and let
$Y_j=\{x\in D_A;\  g^j_\ell(a)x=0, \ell=1,\ldots, \ell_j\}$.  If $x\in X_j$, then
$x=\phi_j(a)x$, and since moreover
$g^j_\ell\phi_j=\sum_k h_k f_k\phi_j$ for some $h_k$, it follows that $x\in Y_j$.
Thus $X_j\subset Y_j$. Furthermore, if for each $j$ the common zero set of $g_\ell^j$ is
just the point $\alpha_j$, then by \eqref{ink}, $\sigma(a,\overline Y_j)\subset\{\alpha_j\}$.
If $x\in Y_j\cap Y_i$, therefore $\sigma_x(a)=\emptyset$, and hence $x=0$ since $a$
is a hyperoperator.  It therefore follows that $X_j=Y_j$.

If all zeros of $f$ are of first order, i.e., the local ideal at $\alpha^j$
is generated by $\xi_i-\alpha^j_i$, $i=1,\ldots,n$, then $X_j$ is the eigenspace
$$
X_j=\{x\in D_A;\  a_ix=\alpha^j_ix, \  i=1,\ldots,n\}.
$$

\smallskip

If $A\in H_{\D(\C)}(X)$ and 
$f$ is holomorphic in a neighborhood of $\sigma(a)$ and the zeros $\alpha^j$ 
have multiplicities $r_j$, then
$
X_j=\{x\in D_A;\   (a-\alpha_j)^{r_j}x=0\}.
$
\end{ex}

The situation in this  example  appears naturally when we consider
homogeneous solutions to an equation like  $f(a)x=0$.

\begin{ex}
Let $\A_h(\R^n)$ be a Beurling algebra, cf., Section~\ref{ultra}, 
containing cut-off functions, 
and let $X$ be the
space of inverse Fourier transforms of the dual space $\A_h'$. Then the tuple of commuting 
operators $a_j=i\partial/\partial\xi_j$ on $X$ admits 
an $\A_h$ functional calculus (since $\A_h$ is an algebra).
Then 
$D_a$ is the space of (inverse) Fourier transforms of elements with compact supports 
in $\A_h'$. Notice that
$\A_h'$ contains all distributions with compact support, but also some hyperfunctions 
of infinite
order.   Let $f$ be a  $\A_h$-smooth  mapping and consider the space
$\{x\in D_a;  \ f(a)=0\}$.  If $x$ is the inverse Fourier transform of $u$, then 
$f(t)u(t)=0$, which
means that $u$ has support on $Z=\{t\in\R^n;\  f(t)=0\}$.   It follows that we have the  
representation
\begin{equation}\label{eren}
x(\xi)=\int_{\R^n} e^{i\xi\cdot  t}u(t)dt,
\end{equation}
meaning the action of  $u$ on $t\mapsto \exp{i\xi\cdot t}$. Since $u$ has support on 
the set $Z=\{ f=0\}$, $x$ is  expressed as a combination of exponentials with frequencies
in $Z$.
\end{ex}

Even if $f$ is a polynomial, only  solutions generated by real frequencies can appear
as long as we have restricted to non-quasianalytic classes. 
To get an operator-theoretic frame of this kind
for the general   fundamental principle of Ehrenpreis and Palamodov,
\cite{Eh} and \cite{Pal}, one must consider operators that only admit
a holomorphic functional calculus.

\section{Non-commuting hyperoperators}

Assume that $A_1,\ldots,A_m$ are in $H_{\D(\R^{n_j})}$ but not necessarily
commuting. Then we can form the tensor product $A=A_1\otimes\cdots \otimes A_m$,
and obtain a linear continuous, though not multiplicative,
operator $\D(\R^n)\to\L(X)$, where $n=n_1+\cdots +n_m$.
This can also be done explicitly by the formula
\eqref{defa}. We also write this
operator of course as $\phi(a_1,\ldots,a_m)$. 
In case when all $n_j=1$ and $a_j$ are  HS operators, we get back the 
definition in \cite{AS}. 
Now the order of the operators is crucial. 
Therefore it is convenient to use  Feynman notation, see, e.g., \cite{NSS}. Then
this operator $\phi(a_1,\ldots,a_m)$ can be written
$\phi(\sr{m}{a_1},\ldots,\sr{1}{a_m})$ indicating that the operator
$a_m$ is to be applied first, then $a_{m-1}$ etc and finally $a_1$,
and the order is reflected by the order of the resolvents. Therefore,
if $b$ is a bounded operator one can easily define for instance
$$
\phi(\sr{3}{a_1},\sr{1}{a_2})\sr{2}{b_{}}=
\int\int(\partial_{\bar {\zeta_1}}\partial_{\bar
{\zeta_2}})
\tilde\phi(\zeta_1,\zeta_2)\w
\omega^{\chi_1}_{\zeta_1-a_1}\w b \omega^{\chi_1}_{\zeta_2-a_2}.
$$
Notice that this is {\it not} an ordinary composition of
$f(\sr{3}{a_1},\sr{1}{a_2})$
and $b$, while for instance
$
f(\sr{2}{a_1},\sr{1}{a_2})\sr{3}{b_{}}=b f(\sr{2}{a_1},\sr{1}{a_2}).
$

\end{document}